\documentclass[10pt]{amsart}
\usepackage{euscript,amssymb,xspace}
\usepackage[author-year]{amsrefs}

   \title[Dynamic Walks]{On dynamical Gaussian random walks
      }
   \thanks{The research of
      D.\@ Kh.\@ is partially supported by a grant from the NSF}
   \address{Department\@ of Mathematics\\The University\@ of Utah\\
      155 S.\@ 1400 E.\\Salt Lake City, UT 84112--0090}
   \author[D. Khoshnevisan]{Davar Khoshnevisan}
   \email{davar@math.utah.edu}
   \urladdr{http://www.math.utah.edu/\~{}davar}
   \author[D.\@ A.\@ Levin]{David A.\@ Levin}
   \address{Department\@ of Mathematics\\The University\@ of Utah\\
      155 S.\@ 1400 E.\\Salt Lake City, UT 84112--0090}
   \email{levin@math.utah.edu}
   \urladdr{http://www.math.utah.edu/\~{}levin}
   \author[P.\@ M\'endez]{Pedro J. M\'endez-Hern\'andez}
   \address{Department\@ of Mathematics\\The University\@ of Utah\\
      155 S.\@ 1400 E.\\Salt Lake City, UT 84112--0090}
   \email{mendez@math.utah.edu}
   \urladdr{http://www.math.utah.edu/\~{}mendez}

\theoremstyle{plain}{
\newtheorem{theorem}{Theorem}[section]}
\theoremstyle{plain}{
}
\theoremstyle{plain}{
   \newtheorem{lemma}[theorem]{Lemma}}
\theoremstyle{plain}{
   \newtheorem{corollary}[theorem]{Corollary}}
\theoremstyle{definition}{
   }
\theoremstyle{definition}{
   }
\theoremstyle{remark}{
   \newtheorem{remark}[theorem]{Remark}}
\numberwithin{equation}{section}
\newcommand{\F}{\mathfrak{F}}
\newcommand{\N}{\mathfrak{N}}

\newcommand{\e}{\varepsilon}
\newcommand{\s}{\sigma}

\renewcommand{\P}{\mathsf{P}}

\newcommand{\ee}{\mathsf{e}}
\newcommand{\PN}{\P_{_\N}}
\newcommand{\E}{\mathsf{E}}

\newcommand{\EN}{\E_{_\N}}

\newcommand{\VN}{\text{\rm Var}_{_\N}}
\newcommand{\R}{\mathbb{R}}

\renewcommand{\l}{\lambda}

\newcommand{\Given}{\; \vline \;}

\newcommand{\al}{\alpha}
\newcommand{\cadlag}{c\`adl\`ag\xspace}
\begin{document}
\begin{abstract}
   Motivated by the recent work of Benjamini, H\"aggstr\"om,
   Peres, and Steif \ycite{benjamini}
   on dynamical random walks, we: (i) Prove that, after a
   suitable normalization,
   the dynamical Gaussian walk converges weakly to the Ornstein--Uhlenbeck
   process in classical Wiener space; (ii) derive sharp tail-asymptotics
   for the probabilities of large deviations of the said dynamical walk;
   and (iii) characterize (by way of an integral test)
   the minimal envelop(es) for the growth-rate of the
   dynamical Gaussian walk. This development also
   implies the tail capacity-estimates of Mountford
   \ycite{mountford} for
   large deviations in classical Wiener space.

   The results of this paper give a partial affirmative
   answer to the problem, raised in \ocite{benjamini}*{Question 4}
   of whether there are precise connections
   between the OU process in classical
   Wiener space and dynamical random walks.
\end{abstract}
\maketitle
\section{Introduction and Main Results}

Let $\{\omega_j\}_{j=1}^\infty$ denote a sequence of i.i.d.\
random variables, and to each $\omega_j$ we associate a
rate-one Poisson process with jump times
$0<\tau_j(1)<\tau_j(2)<\ldots$\,. (All of the said processes
are assumed to be independent from
one another.) Now at every jump-time of the $j$th
Poisson process, we replace the existing $\omega$-value by an
independent copy. In symbols,
let $\{\omega_j^k\}_{j,k=1}^\infty$ be a double-array
of i.i.d.\ copies of the $\omega_j$'s---all independent of
the Poisson clocks---and define the process
$X:=\{X_j(t);\ t\ge 0\}_{j=1}^\infty$ as follows: For all $j\ge 1$,
\begin{equation}\begin{split}
   X_j(0) & := \omega_j,\\
   X_j (t) & :=  \omega_j^k,\quad {}^\forall
   t\in\left[\tau_j(k), \tau_j(k+1)\right).
\end{split}\end{equation}

We remark that, as a process indexed by $t$,
$t\mapsto (X_1(t),X_2(t),\cdots)$ is a stationary Markov
process in $\R^\mathbb{N}$ whose invariant measure is the
product measure $\mu^\infty$, where $\mu$ denotes the law of $\omega_1$.

Recently,
Benjamini, H\"aggstr\"om, Peres, and Steif \ycite{benjamini}
have introduced
\emph{dynamical random walks} as the partial-sum processes
that are associated to the Markov process $X$.
In other words, the dynamical walk associated to the distribution
$\mu$ is defined as the two-parameter process 
$S:=\{S_n(t)\}_{n\ge 1, t\ge 0}$
that is defined by
\begin{equation}
   S_n(t) := X_1(t) +\cdots + X_n(t),
   \qquad{}^\forall n\ge 1,\ t\ge 0.
\end{equation}
{F}rom now on, we specialize our dynamical walks by assuming that
the incremental distribution $\mu$ is standard normal, i.e., for
all $x\in\R$,
\begin{equation}
   \mu\left( [x,\infty)\right) = 1-\Phi(x) :=
   \bar\Phi (x) := \int_x^\infty \frac{e^{-z^2/2}}{\sqrt{2\pi}}
   \, dz.
\end{equation}
Our forthcoming analysis depends on this simplification
in a critical way.

Now consider the following rescaled dynamical Gaussian walk $U^n$:
\begin{equation}
   U^n_t(s):= \frac{S_{\lfloor nt\rfloor}(s)}{\sqrt{n}},
   \qquad{}^\forall s,t\in[0,1]).
\end{equation}
Our first contribution is the following large-sample result
on dynamical Gaussian walks.
\begin{theorem}\label{thm:weak}
   As $n$ tends to infinity, the random field $U^n$ converges weakly in
   $D([0,1]^2)$ to the continuous centered Gaussian random field
   $U$ whose covariance is
   \begin{equation}
      \E\left\{ U_s(t) U_{s'}(t') \right\}= e^{-|s-s'|} \min(t,t'),
      \qquad {}^\forall s,s',t,t'\in[0,1].
   \end{equation}
\end{theorem}
\noindent(For information on $D([0,1]^2)$
consult Section~\ref{sec:WeakConv}.)

Before proceeding further,
we make two tangential remarks.

\begin{remark}
   The limiting
   random field $U$ has the following interpretation:
   \begin{equation}\label{eq:U}
      U_t(s) := e^{-s}B\left( e^{2s},t \right),\qquad
      {}^\forall s,t\in[0,1],
   \end{equation}
   where $B$ is the two-parameter Brownian sheet. Standard arguments
   then show that $\mathcal{U}:=\{ U_t\}_{t\ge 0}$ is
   an infinite-dimensional stationary diffusion on the classical
   Wiener space $C([0,1])$, and the invariant measure of $\mathcal{U}$
   is, in fact, the Wiener measure on $C([0,1])$. The process
   $\mathcal{U}$ is the so-called \emph{Ornstein--Uhlenbeck} (OU) process
   in classical Wiener space. Theorem~\ref{thm:weak},
   in conjunction with this observation, gives a partial affirmative
   answer to \ocite{benjamini}*{Question 4}, where it is asked whether there
   are precise potential-theoretic
   connections between the dynamical (here, Gaussian) walks, and
   the OU process in $C([0,1])$.
\end{remark}

\begin{remark}
   Theorem~\ref{thm:weak} can be viewed as a construction of
   the OU process in $C([0,1])$. This is an interesting
   process in and of itself, and
   arises independently in diverse areas
   in stochastic analysis. For three samples, see 
   \ocite{kuelbs}, \ocite{malliavin}, and \ocite{walsh}.
   The elegant relation (\ref{eq:U}) to the Brownian sheet was
   noted by David Williams;
   cf.\ 
   \ocite{meyer}*{appendix}.
\end{remark}

Our next result elaborates further on the connection between
the dynamical Gaussian walk and the process $\mathcal{U}$.

\begin{theorem}\label{thm:key-est}
   Choose and fix a sequence $\{ z_j\}_{j=1}^\infty$ that
   satisfies
   \begin{equation}\label{eq:z}
      \inf_n z_n \ge 1,\
      \lim_{n\to\infty}z_n=+\infty,\ \text{ and }\
      \lim_{n\to\infty}\sqrt{\frac{\log n}{n}} z_n=0.
   \end{equation}
   Then, as $n\to\infty$,
   \begin{equation}\label{eq:key}
      \frac{1+o(1)}{9} z_n^2 \bar\Phi(z_n)\le
      \P\left\{ \sup_{t\in[0,1]}
      S_n(t) \ge z_n\sqrt{n} \right\}
      \le (2+o(1)) z_n^2 \bar\Phi(z_n).
   \end{equation}
\end{theorem}
The following reformulation of a theorem of Mountford \ycite{mountford}
provides the analogue for the standard OU process $U:=\{ U_1(s);\
s\ge 0\}$: \emph{There exists a constant
$K_{\ref{eq:mountford}}>1$ such that
\begin{equation}\label{eq:mountford}
   K_{\ref{eq:mountford}}^{-1} z^2 \bar\Phi(z) \le
   \P\left\{ \sup_{s\in [0,1]} U_1(s)
   \ge z \right\} \le K_{\ref{eq:mountford}} z^2 \bar\Phi(z),
   \qquad{}^\forall z\ge 1.
\end{equation}
}
{F}or a refinement see \ocite{pickands}, and
also \ocite{qualls1}.

The apparent similarity between Theorem~\ref{thm:key-est} and
(\ref{eq:mountford})
is based on more than mere analogy. Indeed, Theorems~\ref{thm:weak}
and~\ref{thm:key-est} together imply (\ref{eq:mountford}) as a corollary.
This can be
readily checked; cf.\ the last line of \S\ref{subsec:D}.

As a third sample from our present work, we show a pathwise
implication of Theorem~\ref{thm:key-est}. This is
the dynamical analogue of the celebrated ``integral test'' of
Erd\H{o}s \ycite{erdos}. Define
the map $\EuScript{J}(H)$, for all nonnegative measurable
functions $H$, by
\begin{equation}\label{eq:JJ}
   \EuScript{J}(H) := \int_1^\infty \frac{%
   H^4(t) \bar\Phi(H(t))}{t}\, dt.
\end{equation}

\begin{theorem}\label{thm:erdos}
   Suppose that $H$ is a nonnegative nondecreasing function.
   Then:
   \begin{enumerate}
      \item[(i)]
         If $\EuScript{J}(H)<+\infty$,
         then with probability one,
         \begin{equation}
            \sup_{t\in[0,1]} S_n(t) < H(n)\sqrt{n},\
            \text{ for all but a finite number of $n$'s}.
         \end{equation}
      \item[(ii)]
         Conversely, if $\EuScript{J}(H)=+\infty$, then  with
         probability one there exists a $t\in[0,1]$,
         such that
         \begin{equation}
            S_n(t) \ge H(n)\sqrt{n},\
            \text{ for an infinite number of $n$'s}.
         \end{equation}
   \end{enumerate}
\end{theorem}

\begin{remark}
   Owing to (\ref{eq:Phiasymp}) below, we have
   \begin{equation}\label{eq:sumint2}
      \EuScript{J}(H)<+\infty\ \Longleftrightarrow\
      \int_1^\infty H^3(t) e^{-\frac12 H^2(t)}\,
      \frac{dt}{t}<+\infty.
   \end{equation}
\end{remark}

We recall that the Erd\H{o}s integral test
asserts that
$S_n(0)>H(n)\sqrt{n}$ for infinitely many $n$ (a.s.)
if and only if $\int_1^\infty H(t) e^{-\frac12 H^2(t)}
t^{-1}\, dt<+\infty$.
Combining the preceding remark with Theorem~\ref{thm:erdos}
immediately leads us to the following result whose elementary
proof is omitted.

\begin{corollary}
   Given $\tau\in[0,1]$, 
   \begin{equation}
      \limsup_{n\to\infty} \frac{\left[S_n (\tau)\right]^2- 2n \ln\ln n }{%
      n\ln\ln\ln n} = 3,\qquad\text{a.s.}
   \end{equation}
   On the other hand, there exists a (random) $T\in[0,1]$, such
   that
   \begin{equation}
      \limsup_{n\to\infty} \frac{\left[S_n (T)\right]^2- 2n \ln\ln n }{%
      n\ln\ln\ln n} = 5,\qquad\text{a.s.}
   \end{equation}
\end{corollary}

\begin{remark}
   In the terminology of \ocite{benjamini}, our
   Theorem~\ref{thm:erdos} has the consequence that
   the Erd\H{o}s characterization of the upper class
   of a Gaussian random walk is ``dynamically sensitive.''
   This is in contrast
   to the fact that the LIL itself is ``dynamically stable.''
   In plain terms, the latter means that with probability one,
   \begin{equation}\label{eq:LIL}
      \limsup_{n\to\infty} \frac{S_n(t)}{\sqrt{2n\ln\ln n}}
      =1,\qquad{}^\forall t\in[0,1].
   \end{equation}
   See \ocite{benjamini}*{Theorem 1.2}.
\end{remark}

The organization of this paper is as follows:
In \S\ref{sec:clocks} we state and prove a
theorem on the Poisson clocks that, informally speaking,
asserts that with overwhelming
probability the typical clock is at mean-field all the time,
and this happens simultaneously
``over a variety of scales.'' This material may be
of independent technical interest to the reader.

In \S\ref{sec:regression},
we make a few computations with Gaussian random variables.
These calculations are simple consequences of
classical regression analysis of mathematical statistics,
but since we need the exact forms of the ensuing estimates, we include some
of the details.

After a brief discussion of the space $D([0,1]^2)$,
Theorem~\ref{thm:weak} is then proved in \S\ref{sec:WeakConv}.
Our proof relies heavily on the general machinery of
Bickel and Wichura \ycite{bickel}.

Theorem~\ref{thm:key-est} is more difficult to prove;
its proof is split across \S\ref{sec:quenched-UB},
\S\ref{sec:quenched-LB}, and \S\ref{sec:key}. The key idea
here is that estimates, similar to those in Theorem~\ref{thm:key-est},
hold in the quenched setting, where the implied
conditioning is made with respect to the clocks.

{F}inally, we derive Theorem~\ref{thm:erdos}
in \S\ref{sec:erdos}.
Our proof combines Theorem~\ref{thm:key-est},
a localization trick,
and the combinatorial method of \ocite{erdos}.

Throughout,
we frequently use the elementary facts that for all $y>0$,
\begin{equation}\label{eq:Phiasymp}
   \bar\Phi(y)\le e^{-y^2/2},\  \text{and}\
   \bar\Phi(z)=\frac{1+o(1)}{z\sqrt{2\pi}} e^{-z^2/2}\qquad
   (z\to\infty).
\end{equation}
We have used Bachmann's
``little-$o$/big-$O$'' notation to simplify
the exposition.

\medskip
\noindent\textbf{Acknowledgment.} We are grateful to
Yuval Peres for introducing us to this subject, and for
a number of interesting discussions.

\section{Regularity of the Clocks}
\label{sec:clocks}

Consider the random field
$\{ N^n_{s\to t};\
0\le s\le t,\, n\ge 1\}$ that is defined as follows: Given $s\le t$
and $n\ge 1$,
$N^n_{s\to t}$ denotes the Poisson-based number of changes made from
time $s$ to time $t$; i.e.,
\begin{equation}
   N^n_{s\to t} := \sum_{j=1}^n \mathbf{1}_{\{X_j(t)\neq X_j(s)\}}.
\end{equation}

It is clear that $N^n_{s\to t}$ is a sum of $n$ i.i.d.\
$\{0,1\}$-valued random variables. Because
we know also that $\P\{ X_1(s)=X_1(t)\}=e^{-|t-s|}$, we
can deduce from the strong law for such binomials that for $n$ large,
$N^n_{s\to t} \simeq n (1-e^{-|t-s|})$. The following is an
estimate that ensures that,
in the mentioned approximation, a good amount of
uniformity in $s$ and $t$ is preserved.

\begin{theorem}\label{thm:LD}
   If $\{ \Delta_j\}_{j=1}^\infty$ is a sequence in $[0,1]$
   such that $\lim_{n\to\infty}\Delta_n=0$, then
   for all $n\ge 1$ and $\alpha\in(0,1)$,
   \begin{equation}
      \P\left\{ \sup_{\substack{0\le s\le t \le 1: \\
       t -s\ge \Delta_n}}
      \left| \frac{N^n_{s\to t}}{%
      \E  N^n_{s\to t} } - 1
      \right| \ge \alpha \right\}\leq \frac{512
      }{\alpha^2 \Delta_n^2}
      \exp\left(- \frac{3\alpha^3 n \Delta_n}{2304}\right),
   \end{equation}
   where $\sup\varnothing:=0$.
\end{theorem}

This, and the Borel--Cantelli lemma, together imply the following
result that we shall need later on. In rough terms,
it states that as long as the ``window size''
is not too small, then the Poisson clocks are mean-field.
\begin{corollary}\label{cor:LD}
   If $\Delta_n\to 0$ in $[0,1]$ satisfies
   $\lim_{n\to\infty}n(\log n)^{-1}\Delta_n=+\infty$, then
   with probability one,
   \begin{equation}
      \lim_{n\to\infty} \sup_{\substack{0\le s\le t \le 1:\\
       t-s\ge \Delta_n}}
      \left| \frac{N^n_{s\to t}}{\E N^n_{s\to t}
      } - 1 \right| =0.
   \end{equation}
\end{corollary}
It is not hard to convince oneself that the preceding fails
if the ``window size'' $\Delta_n$ decays too rapidly.

\begin{proof}[Proof of Theorem~\ref{thm:LD}]
   Throughout this proof, $\alpha\in(0,1)$ is held fixed.

   We first try to explain the significance of the condition
   $t-s\ge\Delta_n$ by obtaining a simple lower bound on
   $\E N^n_{s\to t}$ in this case.

   Observe the following simple bound:
   \begin{equation}\label{eq:expbound}
      \frac{x}{2} \le 1- e^{-x} \le x,\qquad{}^\forall x\in[0,1].
   \end{equation}
   This shows that
   \begin{equation}\label{eq:ENLower}
      \inf_{\substack{0\le s\le t\le 1:\\
      t-s\ge\Delta_n}}
      \E N^n_{s\to t} \ge
      \frac{n\Delta_n}{2}.
   \end{equation}

   Next we recall an elementary large deviations bound for
   Binomials.
   According to Bernstein's inequality (cf.\ 
   \ocite{bennett};
   also see the elegant inequalities of
   Hoeffding \ycite{hoeffding}),
   if $\{ B_j\}_{j=1}^\infty$ are i.i.d.\
   Bernoulli random variables with $\P\{ B_1=1\}:=p$, then
   \begin{equation}\label{eq:bernstein}
      \P\left\{ \left| B_1+\cdots+B_n -np \right|\ge n\lambda\right\}
      \le 2\exp\left( - \frac{n\lambda^2}{
      2p+ \frac23 \lambda} \right).
   \end{equation}
   Apply this with $B_j:=\mathbf{1}_{\{ X_j(s)\neq X_j(t)\}}$,
   for arbitrary $s\le t$ and $\lambda:=\alpha[ 1-e^{-(t-s)}]$,
   to deduce that for all $\alpha\in(0,1)$ and $n\ge 1$,
   \begin{equation}\begin{split}
      &\P\left\{ \left| N^n_{s\to t} - \E
         N^n_{s\to t} \right| \ge\alpha
         \E N^n_{s\to t} \right\}\\
      &\quad\le 2 \exp\left( -\frac{ \alpha^2 n
         \left[ 1-e^{-(t-s)}\right]}{%
         2 +\frac23 \alpha}\right)\\
      &\quad\le 2 \exp\left( -\frac{ 3\alpha^2 n
         \left[ 1-e^{-(t-s)}\right]}{%
         8}\right).
   \end{split}\end{equation}
   From (\ref{eq:expbound}) we can deduce that for all $\alpha\in(0,1)$
   and $n\ge 1$,
   \begin{equation}\label{eq:unifclock1}
      \sup_{\substack{ 0\le s\le t\le 1:\\
      |s-t|\ge\Delta_n}}
      \P\left\{ \left| N^n_{s\to t} - \E
      N^n_{s\to t} \right| \ge\alpha
      \E N^n_{s\to t} \right\}
      \le 2 \exp\left( - \frac{3\alpha^2 n\Delta_n}{16}\right).
   \end{equation}

   Next, we choose and fix integers $k_1<k_2<\cdots\to\infty$
   as follows:
   \begin{equation}\label{eq:k_nDelta_n}
      k_n := \left\lfloor 1+\frac{8}{\alpha \Delta_n}
      \right\rfloor\quad\text{ so that}\quad  
      \frac{\alpha\Delta_n}{9} \le k_n^{-1} \le
      \frac{\alpha\Delta_n}{8}.
   \end{equation}
   Based on these, we define
   \begin{equation}
      \Gamma_n := \left\{ \frac{j}{k_n};\ 0\le j\le k_n\right\}.
   \end{equation}
   Then it follows immediately from (\ref{eq:unifclock1}) 
   and (\ref{eq:k_nDelta_n}) that
   \begin{equation}\label{eq:Fest}
      \P\left\{ \sup_{\substack{ 0\le s\le t\le 1: \\
      s,t\in \Gamma_n}}
      \left| \frac{N^n_{s\to t} }{\E N^n_{s\to t}}
      -1\right| \ge\alpha \right\}
      \le (k_n+1)^2 \exp\left(- 
      \frac{3\alpha^3 n\Delta_n}{144}\right).
   \end{equation}

   Given any point $u\in[0,1]$, define
   \begin{equation}\begin{split}
      \underline{u}_n & := \max\left\{ r\in[0,u]:\
         r\in \Gamma_n \right\}\\\
      \overline{u}_n & := \min\left\{ r\in[u,1]:\
         r\in \Gamma_n\right\}.
   \end{split}\end{equation}
   These are the closest points to $u$ in $\Gamma_n$ from
   below and above respectively.
   We note, in passing,
   that $0\le \overline{u}_n - \underline{u}_n\le k_n^{-1}$.
   Moreover, thanks to (\ref{eq:k_nDelta_n}), whenever
   $0\le s\le t\le 1$ satisfy $t-s\ge\Delta_n$, it follows that
   $\overline{s}_n < \underline{t}_n$ with room to spare.
   We will use this fact
   without further mention. Moreover, for such a pair $(s,t)$,
   \begin{equation}\label{eq:3N}
      N^n_{\overline{s}_n\to \underline{t}_n}  \le
      N^n_{s\to t} \le
      N^n_{\underline{s}_n\to \overline{t}_n}.
   \end{equation}
   This follows from the fact that
   with $\P$-probability one,
   once one of the $X_j(u)$'s is updated, then from that point
   on it will never be
   replaced back to its original state. (This is so because the chances
   are zero that two independent normal
   variates are equal to one another.) The preceding display motivates
   the following bound: For all $0\le s\le t\le 1$,
   \begin{equation}\begin{split}
      \E\left\{ \left| N^n_{\underline{s}_n \to \overline{t}_n} -
         N^n_{\overline{s}_n\to \underline{t}_n}\right| \right\}
         & = n e^{-(\underline{t}_n - \overline{s}_n)}
         \left[ 1 - e^{-(\overline{t}_n-\underline{t}_n) -
         (\overline{s}_n-\underline{s}_n)}\right]\\
      & \le \frac{2n}{k_n},
   \end{split}\end{equation}
   where the last inequality follows from (\ref{eq:expbound}).
   Owing to (\ref{eq:ENLower}) and (\ref{eq:k_nDelta_n}), we have
   the crucial estimate,
   \begin{equation}\label{eq:ENControl}
      \sup_{\substack{ 0\le s\le t\le 1:\\
      t-s\ge\Delta_n} }
      \E\left\{ \left| N^n_{\underline{s}_n \to \overline{t}_n} -
      N^n_{\overline{s}_n\to \underline{t}_n}\right| \right\}
      \le \frac{\alpha}{2}  \inf_{\substack{ 0\le u\le v\le 1:\\
      v-u\ge\Delta_n}}
      \E N^n_{u\to v}.
   \end{equation}
   This and (\ref{eq:3N}) together imply the following bound
   uniformly for all $0\le s\le t\le 1$ that satisfy
   $t-s\ge\Delta_n$:
   \begin{equation}
      \left| \widetilde{N}^n_{s\to t} \right|
      \le \frac{\alpha}{2} \inf_{\substack{ 0\le u\le v\le 1 \\
      v-u\ge\Delta_n}}
      \E N^n_{u\to v} + \max\left(
      \left|\widetilde{N}^n_{%
      \overline{s}_n\to \underline{t}_n}\right| ~,~
      \left|\widetilde{N}^n_{\underline{s}_n \to \overline{t}_n}
      \right|\right),
   \end{equation}
   where $\widetilde{Z}:=Z-\E Z$ for any integrable
   random variable $Z$. Therefore,
   \begin{equation}\begin{split}
      & \P\left\{ {}^\exists t-s\ge\Delta_n:\ \left|
         \widetilde{N}^n_{s\to t} \right| \ge\alpha
         \E N^n_{s\to t} \right\}\\
      &  \le \P\left\{ {}^\exists t-s\ge\Delta_n:\ \max\left(
         \left|\widetilde{N}^n_{\overline{s}_n\to
         \underline{t}_n}\right| ~,~
         \left|\widetilde{N}^n_{\underline{s}_n \to \overline{t}_n}
         \right|\right)\ge \frac{\alpha}{2}
         \E N^n_{s\to t} \right\}.
   \end{split}\end{equation}
   Another application of (\ref{eq:ENControl}) yields
   \begin{equation}\begin{split}
      & \P\left\{ {}^\exists t-s\ge\Delta_n:\ \left|
         \widetilde{N}^n_{s\to t} \right| \ge\alpha
         \E N^n_{s\to t} \right\}\\
      &  \le \P\left\{ {}^\exists t-s\ge\Delta_n:\
         \left|\widetilde{N}^n_{\overline{s}_n\to \underline{t}_n}\right|
         \ge \frac{\alpha}{2}\left( 1 -\frac{\alpha}{2}
         \right)
         \E N^n_{\overline{s}_n\to \underline{t}_n}
         \right\}\\
      &\quad + \P\left\{ {}^\exists t-s\ge\Delta_n:\
         \left|\widetilde{N}^n_{\underline{s}_n
         \to \overline{t}_n}\right|
         \ge \frac{\alpha}{2}\left( 1 -\frac{\alpha}{2} \right)
         \E N^n_{\underline{s}_n \to
         \overline{t}_n} \right\}\\
      & \le 2 \P\left\{ \max_{\substack{ 0\le u\le v\le 1:\\
         u,v\in \Gamma_n}} \left|
         \frac{N^n_{u\to v}}{\E  N^n_{u\to v}}- 1\right|
         \ge \frac{\alpha}{4} \right\}\\
      & \le 2
         \left(k_n^2 +1\right)
         \exp\left(- \frac{3\alpha^3 n\Delta_n}{2304}\right),
   \end{split}\end{equation}
   owing to (\ref{eq:Fest}). Because
   $k_n+1 \le 16(\alpha\Delta_n)^{-1}$, this proves the theorem.
\end{proof}

\section{A Little Regression Analysis}
\label{sec:regression}

Define $\F_t^n$ to be the augmented right-continuous
$\s$-algebra generated by the variables $\{ S_n(v);\, v\le t\}$
and $\N$, where the latter is the $\s$-algebra
generated by all of the Poisson clocks. For convenience,
we write $\PN\{\cdots\}$ and $\EN\{\cdots\}$ in
place of $\P\{\cdots\,|\,\N\}$ and $\E\{\cdots\,|\,\N\}$,
respectively. We refer to $\PN$ as a random ``quenched'' measure,
and $\EN$ is its corresponding expectation operator.
We will also write $\VN$ for the corresponding conditional
variance.

\begin{lemma}\label{lem:covariance} If $0\le u\le v$, then
   the following hold $\P$-almost surely: For all $x\in\R$,
   \begin{equation}\begin{split}
      \EN\left\{ S_n(v) \, \Big|\, S_n(u)=x \right\} & =
         \left( 1 - \frac{N^n_{u\to v}}{n}\right)x,\\
      \VN \left( S_n(v)\, \Big|\,
         S_n(u) =x \right) & = N^n_{u\to v} \left[ 2 - \frac{
         N^n_{u\to v}}{n}  \right].
   \end{split}\end{equation}
\end{lemma}

\begin{proof}
   From time $u$ to time $v$, $N^n_{u\to v}$-many
   of the increments are changed; the remaining
   $(n-N^n_{u\to v})$ increments are left unchanged. Therefore,
   we can write
   \begin{equation}\label{eq:decompose}\begin{split}
      S_n(u) & = V_1 + V_2\\
      S_n(v) & = V_1 + V_3,
   \end{split}\end{equation}
   where: (i) $V_1$, $V_2$, and $V_3$ are independent;
   (ii) the distribution of $V_1$ is the same as that
   of $S_{n-N^n_{u\to v}}(0)$; and (iii)
   $V_2$ and $V_3$ are identically distributed and
   their common distribution is that of
   $S_{N^n_{u\to v}}(0)$. The result follows
   from standard calculations
   from classical regression analysis.
\end{proof}

This immediately yields the following.

\begin{lemma}\label{lem:cond}
   For all $x,y\ge 0$, all times $0\le u\le v$, and all integers $n\ge 1$,
   \begin{equation}\begin{split}
      \PN \left\{ S_n(v) \ge y \,\Big|\, \F_u^n \right\}
         &= \PN \left\{ S_n(v) \ge y \,\Big|\, S_n(u) \right\} \\
      & = \bar\Phi\left(
         \frac{y -  \left( 1-\frac1n N^n_{u\to v} \right)S_n(u)}{
         \sqrt{N^n_{u\to v}\left( 2- \frac1n N^n_{u\to v}\right)}}
         \right),\qquad\P\text{-a.s.}
   \end{split}\end{equation}
\end{lemma}

We will also have need for
the following whose elementary proof we omit.

\begin{lemma}\label{lem:phi-bar}
   For all $z\ge 1$ and $\e>0$, we have
   $\bar\Phi( z + \e z ) \le  e^{- z^2\e} \bar\Phi(z).$
\end{lemma}

Next is a ``converse'' inequality.
Unlike the latter lemma, however, this one merits a brief
derivation.

\begin{lemma}\label{eq:phi-bari}
   If $\gamma>0$, then
   \begin{equation}
      \bar\Phi\left( z -\frac{\gamma}{z} \right)
      \le \left( 1 + e^{2\gamma}\right)
      \bar\Phi(z), \qquad {}^\forall z \ge\sqrt{\gamma}.
   \end{equation}
\end{lemma}

\begin{proof}
   We make a direct computation:
   \begin{equation}\begin{split}
       \bar\Phi\left( z -\frac{\gamma}{z} \right) & =
          \frac{1}{\sqrt{2\pi}}\int_z^\infty \exp\left\{
          -\frac12 \left( y -\frac{\gamma}{z} \right)^2\right\}\, dy\\
       &\le \frac{1}{\sqrt{2\pi}} \int_z^{2z} e^{-\frac12 y^2 }
          e^{\gamma y/z}\, dy + \bar\Phi\left( 2z - 
          \frac{\gamma}{z}\right)\\
       &\le e^{2\gamma}\bar\Phi(z) + \bar\Phi\left( 2z 
          - \frac{\gamma}{z}\right).
   \end{split}\end{equation}
   On the other hand, if $z\ge\gamma /z$, then $2z-\gamma/z\ge z$,
   and so $\bar\Phi(2z-\gamma z^{-1})\le\bar\Phi(z)$.
   This completes the proof.
\end{proof}

\section{Weak Convergence}
\label{sec:WeakConv}

\subsection{The Space $D([0,1]^2)$}
\label{subsec:D}

Let us first recall some facts about the
Skorohod space $D([0,1]^2)$ which was introduced and studied
in \ocite{neuhaus}, \ocite{straf}, and \ocite{bickel}.
Bass and Pyke \ycite{basspyke} provide a theory of weak convergence
in $D(A)$ which subsumes that in $D([0,1]^2)$.

In a nutshell, $D([0,1]^2)$ is the collection of all bounded
functions $f:[0,1]^2\to\R$ such that $f$ is \cadlag with respect
to the partial order $\prec$, where
\begin{equation}
   (s,t)\prec (s',t')\ \Longleftrightarrow\
   s\le s',\ \text{and}\ t\le t'.
\end{equation}
Of course, $f$ is \cadlag with respect to $\prec$ if and only if:
(i) As $(s,t)\downarrow(u,v)$ (with respect to $\prec$),
$f(s,t)\to f(u,v)$; and (ii) if $(s,t)\uparrow(u,v)$, then
$f((u,v)^-):=\lim f(s,t)$ exists.

Once it is endowed with a Skorohod-type metric, the space
$D([0,1]^2)$ becomes a complete separable metric
space \cite{bickel}*{p.\ 1662}.

If $X,X_1,X_2,\ldots$ are random elements of $D([0,1]^2)$,
then $X_n$ is said to converge weakly to $X$
(written $X_n\Rightarrow X$) if for all
bounded continuous functions $\phi:D([0,1]^2)\to\R$,
$\lim_{n\to\infty}\E[\phi(X_n)] =\E[\phi(X)]$. Since the
identity map from $C([0,1]^2)$ onto itself is a topological
embedding of $C([0,1]^2)$ in $D([0,1]^2)$,
if $\phi$ is a continuous functional on $C([0,1]^2)$, then
it is also a continuous functional on $D([0,1]^2)$.

An important example of such a continuous functional is
\begin{equation}
   \phi(x):=\sup_{t\in[0,1]} x(t),\qquad
   {}^\forall x\in D([0,1]^2).
\end{equation}
This example should
provide ample details for deriving
Mountford's theorem (\ref{eq:mountford}) from
Theorems~\ref{thm:weak} and~\ref{thm:key-est}
of the present article.

\subsection{Proof of Theorem~\ref{thm:weak}}
The proof, as is usual in weak convergence,
involves two parts. First, we prove
the convergence of all finite-dimensional distributions. This
portion is done in the quenched setting, for then all processes
involved are Gaussian and we need to compute a covariance or two only.
The more interesting portion is the second part and amounts to
proving tightness. Here we use, in a crucial way,
a theorem of \ocite{bickel}.

\begin{proof}[Proof of Theorem~\ref{thm:weak}] (Finite-Dimensional
   Distributions)
   Given any four (fixed) values of $s,t,s',t'\in[0,1]$,
   \begin{equation}\begin{split}
       \EN\left\{ U^n_t(s) U^n_{t'}(s') \right\} & = \frac1n \EN\left\{
          S_{\lfloor nt\rfloor}(s)
          S_{\lfloor nt'\rfloor}(s') \right\}\\
       & = \frac1n \EN\left\{
          S_{\lfloor nt\rfloor \wedge \lfloor nt'\rfloor}(s)
          S_{\lfloor nt\rfloor \wedge \lfloor nt'\rfloor}(s') \right\}
          \,.
   \end{split}\end{equation}
   Thanks to Lemma~\ref{lem:covariance}, $\P$-almost surely,
   \begin{equation}\begin{split}
       &\EN\left\{ U^n_t(s) U^n_{t'}(s') \right\}\\
       &\quad = \frac1n \left( 1-
          \frac{N^{\lfloor nt\rfloor \wedge \lfloor nt'\rfloor}_{
          (s\wedge s')\to(s\vee s')}}{%
          \lfloor nt\rfloor \wedge \lfloor nt'\rfloor}\right)
          \left( \lfloor nt\rfloor \wedge \lfloor nt'\rfloor \right).
   \end{split}\end{equation}
   On the other hand, by the strong law of large numbers, as $n\to\infty$,
   \begin{equation} \begin{split}
      \frac{N^{\lfloor nt\rfloor \wedge \lfloor nt'\rfloor}_{
         (s\wedge s')\to(s\vee s')}}{ \lfloor nt\rfloor \wedge \lfloor nt'\rfloor
         } & = (1+o(1)) \frac{
         \E N^{\lfloor nt\rfloor \wedge \lfloor nt'\rfloor}_{
         (s\wedge s')\to(s\vee s')}}{%
         \lfloor nt\rfloor \wedge \lfloor nt'\rfloor}\\
      & \to 1 - e^{-|s'-s|},\qquad\text{a.s. [$\P$]}.
   \end{split}\end{equation}
   Therefore, $\P$-almost surely,
   $\lim_{n\to\infty} \EN\{ U^n_t(s) U^n_{t'}(s') \}
   = \E\{ U_t(s) U_{t'}(s') \}$.
   This readily implies that $\P$-almost surely,
   the finite-dimensional distributions of $U^n$ converge
   weakly $[\PN]$ to those of $U$. By the dominated convergence
   theorem, this implies the weak convergence, under $\P$,
   of the finite-dimensional distributions of $U^n$ to those of $U$.
\end{proof}

In order to prove tightness, we appeal to
a refinement to the Bickel--Wichura Theorem 3;
cf.\ \ocite{bickel}*{p.\ 1665}. To do so, we need to
first recall some
of the notation of \ocite{bickel}.

A \emph{block} is a two-dimensional half-open rectangle
whose sides are parallel to the axes; i.e., $I$ is a block
if and only if it has the form $(s,t]\times(u,v]\subseteq( 0,1]^2$.
Two blocks $I$ and $I'$ are \emph{neighboring} if either:
(i) $I=(s,t]\times(u,v]$ and $I'=(s',t']\times(u,v]$ (horizontal
neighboring); or
(ii) $I=(s,t]\times(u,v]$ and $I'=(s,t]\times(u',v']$
(vertical neighboring).

Given any two-parameter stochastic process $Y:=\{ Y(s,t);\
s,t\in[0,1]\}$, and any block $I:=(s,t]\times(u,v]$,
the \emph{increment of $Y$ over $I$} [written as
$\mathcal{Y}(I)$] is defined as
\begin{equation}
   \mathcal{Y}(I) := Y(t,v)-Y(t,u)-Y(s,v)+Y(s,u).
\end{equation}

We are ready to recall the following important result
of~\ocite{bickel}. We have stated it in a way that best suits our
later needs.

\begin{lemma}[Refinement to \ocite{bickel}*{Theorem 3}]
  \label{lem:BW}
   Denote by $\{ Y_n\}_{n\ge 1}$ a sequence of random fields
   in $D([0,1]^2)$ such
   that for all $n\ge 1$, $Y_n(s,t)=0$ if $st=0$. Suppose
   that there exist constants $K_{\ref{lem:BW}}>1$,
   $\theta_1,\theta_2,\gamma_1,\gamma_2>0$ such that
   they are all independent of $n$, and
   whenever $I:=(s,t]\times(u,v]$ and
   $J:=(s',t']\times(u',v']$ are neighboring blocks, and if
   $s,t,s',t'\in n^{-1}\mathbb{Z}\cap[0,1]$, then
   \begin{equation}\label{eq:BW}
      \E\left\{ \left| \mathcal{Y}_n (I)
      \right|^{\theta_1}
      \left| \mathcal{Y}_n(J)  \right|^{\theta_2} \right\} \le
      K_{\ref{lem:BW}} \left| I\right|^{\gamma_1}
      \left| J\right|^{\gamma_2},
   \end{equation}
   where $|I|$ and $|J|$ denote respectively the planar Lebesgue measures
   of $I$ and $J$.
   If, in addition, $\gamma_1+\gamma_2>1$, then $\{ Y_n\}_{%
   n\ge 1}$ is a tight sequence.
\end{lemma}

This is the motivation behind our next lemma which is the
second, and final, step in the proof of Theorem~\ref{thm:weak}.

\begin{lemma}\label{lem:tightness}
   The process $Y_n(t,s):= U^n_t(s)$ satisfies (\ref{eq:BW})
   with the values $K_{\ref{lem:BW}}:= 10$,
   $\theta_1=\theta_2=2$, and $\gamma_1=
   \gamma_2=1$. In particular,
   $\{ U^n\}_{n\ge 1}$ is a tight sequence in $D([0,1]^2)$.
\end{lemma}

\begin{proof}
   We begin by proving that (\ref{eq:BW}) indeed holds with
   the stated constants. This is a laborious, but otherwise
   uninspiring, computation which we include for the sake
   of completeness. This computation is divided into two
   successive steps, one for each possible configuration
   of the neighboring blocks $I$ and $J$.

   \emph{Step 1.}\ (Horizontal Neighboring)
   By stationarity, it suffices to consider only the
   case $I := (0,s]\times(0,u]$
   and $J := (s,t]\times (0,u]$ where $s,t\in n^{-1}\mathbb{Z}$.
   In this case,
   \begin{equation}\begin{split}
      \mathcal{Y}_n (I) &= \frac{S_{ns}(u) -S_{ns}(0)}{\sqrt{n}},\\
      \mathcal{Y}_n (J) &= \frac{S_{nt}(u) - S_{nt}(0)
         -S_{ns}(u) +S_{ns}(0)}{\sqrt{n}},
   \end{split}\end{equation}
   which implies the independence of the two [under $\PN$ and/or $\P$],
   since $k\mapsto S_k$ is a random walk on $D([0,1])$.
   Now, with $\P$-probability one,
   \begin{equation}\label{eq:YI0}
      \EN \left\{ \left| \mathcal{Y}_n(I) \right|^2 \right\}  =
      \frac{2ns- 2\EN\left\{ S_{ns}(u) S_{ns}(0) \right\}}{n}
      = \frac{2 N^{ns}_{0\to u}}{n}.
   \end{equation}
   See Lemma~\ref{lem:covariance}. Therefore,
   $\E \{ | \mathcal{Y}_n(I) |^2 \}
   = 2s [ 1- e^{-u} ]\le 2s u = 2 |I|$.
   By this and the stationarity of the infinite-dimensional random walk
   $k\mapsto S_k$, $\E \{ | \mathcal{Y}_n(J) |^2 \} \le 2|J|$.
   In summary, in this first case of Step 1, we have shown that
   $\E\{|\mathcal{Y}_n(I) \mathcal{Y}_n(J)|^2\}
   \le 4|I|\times |J|$, which is
   certainly less than $10 |I|\times|J|$.

   \emph{Step 2.}\ (Vertical Neighboring)
   By stationarity, we need to consider only the case
   where $I=(0,s]\times(0,u]$ and $J=(0,s]\times(u,v]$,
   where $s\in n^{-1}\mathbb{Z}$. In this case,
   \begin{equation}
      \mathcal{Y}_n (I)  = \frac{S_{ns}(u)-S_{ns}(0)}{\sqrt{n}},\
      \text{and}\
      \mathcal{Y}_n (J)  = \frac{S_{ns}(v)-S_{ns}(u)}{\sqrt{n}}.
   \end{equation}
   These are not independent random variables, and consequently
   the calculations are slightly lengthier in this case.

   Using
   the Markov property and Lemma~\ref{lem:covariance}, we
   $\P$-almost surely have the following:
   \begin{align}
      &\EN \left\{ \left. \left| \mathcal{Y}_n(J) \right|^2 \, \right|\,
         \F^n_u \right\}\nonumber\\
      &\quad = \VN \left( \left. \frac{S_{ns}(v)}{\sqrt{n}}
         \,\right|\, S_{ns}(u)\right) + \left[
         \EN\left\{ \left.
         \frac{S_{ns}(v)-S_{ns}(u)}{\sqrt{n}} \,
         \right|\, S_{ns}(u) \right\} \right]^2 \nonumber \\
      \begin{split}
      & \quad = \frac{N^{ns}_{u\to v}}{n} \left(
         2 - \frac{N^{ns}_{u\to v}}{ns}\right)
         + \left( \frac{N^{ns}_{u\to v}}{ns}
         \right)^2 \, \frac{ \left[ S_{ns}(u)
         \right]^2}{n}\\
      &\quad \le \frac{N^{ns}_{u\to v}}{n}\left[ 2+
         \frac{ \left[ S_{n}(u) \right]^2}{ns}\right].
   \end{split}\end{align}
   In particular, $\P$-almost surely,
   \begin{equation}\label{eq:YIYJ0}\begin{split}
      \EN\left\{ \left| \mathcal{Y}_n(I)
         \right|^2 \left| \mathcal{Y}_n(J) \right|^2
         \right\}
      & = \EN\left\{ \left| \mathcal{Y}_n(I) \right|^2
          \EN\left\{ \left| \mathcal{Y}_n(J) \right|^2
          \Given \F^N_u \right\} \right\}
           \\
      & \le \frac{N^{ns}_{u\to v}}{n} \EN \left\{
        \left| \mathcal{Y}_n(I) \right|^2
        \left[ 2+
         \frac{ \left[ S_{ns}(u) \right]^2}{ns}\right]\right\}\\
      & = \frac{N^{ns}_{u\to v}}{n} \left[
         \frac{4 N^{ns}_{0\to u}}{n}+ \EN \left\{
         \left| \mathcal{Y}_n(I) \right|^2
         \frac{ \left[ S_{ns}(u) \right]^2}{ns}\right\}\right].
   \end{split}\end{equation}
   See (\ref{eq:YI0}) for the last line. Applying the
   Cauchy--Bunyakovsky--Schwarz
   inequality, we obtain
   \begin{equation} \label{eq:furs} \begin{split}
      \EN \left\{
         \left| \mathcal{Y}_n(I) \right|^2
         \frac{ \left[ S_{ns}(u) \right]^2}{ns}\right\}
      & \le \sqrt{\EN \left| \mathcal{Y}_n(I) \right|^4
         \times \EN\left\{
         \frac{ \left[ S_{ns}(u) \right]^4}{n^2s^2}\right\} }\\
      & = \sqrt{3 \EN
         \left| \mathcal{Y}_n(I) \right|^4 },
   \end{split}
   \end{equation}
   since whenever $G$ is a centered Gaussian variate,
   $\E G^4 = 3(\E G^2 )^2$. By applying this identity once
   more in conjunction with (\ref{eq:YI0}), we have
   \begin{equation} \label{eq:Y_again}
     3 \EN \left| \mathcal{Y}_n(I) \right|^4
     \le 9 \left[ \EN \left| \mathcal{Y}_n(I) \right|^2 \right]^2
      = 36 \left[ \frac{N^{ns}_{0 \rightarrow u}}{n} \right]^2 \,.
   \end{equation}
   Plugging \eqref{eq:Y_again} into \eqref{eq:furs} yields
   the following $\P$-almost sure inequality:
   \begin{equation}
      \EN \left\{
      \left| \mathcal{Y}_n(I) \right|^2
      \frac{ \left[ S_{ns}(u) \right]^2}{ns}\right\}
      \le 6 \frac{N^{ns}_{0\to u}}{n}.
   \end{equation}
   We can plug this into (\ref{eq:YIYJ0}) to deduce that $\P$-a.s.,
   \begin{equation}
      \EN\left\{ \left| \mathcal{Y}_n(I)
      \right|^2 \left| \mathcal{Y}_n(J) \right|^2
      \right\} \le 10
      \frac{N^{ns}_{u\to v}}{n} \frac{N^{ns}_{0\to u}}{n} \,.
   \end{equation}
   On the other hand, $N^{ns}_{u\to v}$ and
   $N^{ns}_{0\to u}$ are independent. Therefore,
   \begin{equation}\begin{split}
      \E \left\{ \left| \mathcal{Y}_n(I)
         \right|^2 \left| \mathcal{Y}_n(J) \right|^2
         \right\}
         & \le 10
         \E\left[ \frac{N^{ns}_{u\to v}}{n}\right]
         \E\left[ \frac{N^{ns}_{0\to u}}{n}\right] \\
      &  = 10 s^2
         \left[ 1-e^{-(v-u)}\right]
         \left[ 1-e^{-u}\right] \\
      &  \le 10 su \times s(v-u) \\
      & = 10 |I| \times|J|.
   \end{split}\end{equation}

   We have verified (\ref{eq:BW}) with $K_{\ref{eq:BW}}= 10$,
   $\theta_1=\theta_2=2$, $\gamma_1=\gamma_2=1$.
   Now if it were the case that $Y_n(s,t)=0$ whenever $st=0$,
   we would be done. However, this is not so. To get around this
   small difficulty, note that what we have shown thus far reveals that
   the random fields $(s,t)\mapsto Y_n(s,t)-n^{1/2}S_{ns}(0)$
   ($n=1,2,\ldots$) are tight. On the other
   hand, by Donsker's invariance principle, the processes
   $s\mapsto n^{-1/2}S_{ns}(0)$ ($n=1,2,\ldots$)
   are tight, and the lemma follows from this and the triangle inequality.
\end{proof}
\section{A Quenched Upper Bound}
\label{sec:quenched-UB}

Without further ado,
next is the main result of this section. Note that
it gives quenched tail estimates for
$\sup_{t\in[r,r+1]}S_n(t)$ since the latter
has the same distribution as
$\sup_{t\in[0,1]}S_n(t)$.

\begin{theorem}\label{thm:key-UB}
   Suppose $\{z_j\}_{j=1}^\infty$ is a nonrandom sequence
   that satisfies property (\ref{eq:z}).
   Then with
   $\P$-probability one, for all $ \e>0$,
   there exists an integer $n_0\ge 1$
   such that for all $n\ge n_0$,
   \begin{equation}
      \PN\left\{ \sup_{t\in[0,1]}
      S_n(t) \ge z_n\sqrt{n} \right\}
      \le (2+\e) z_n^2 \bar\Phi(z_n).
   \end{equation}
\end{theorem}

In the remainder of this section we prove Theorem~\ref{thm:key-UB}.
Throughout, we choose and fix a sequence $z_n$ that satisfies
(\ref{eq:z}). Based on these $z_n$'s, we define the ``window size,''
\begin{equation}\label{eq:Delta_n}
   \Delta_n := \frac{1}{16 z_n^2},\qquad{}^\forall n\ge 1.
\end{equation}
 According to (\ref{eq:z}), the
sequence $\{ \Delta_j\}_{j=1}^\infty$ satisfies the conditions of
Theorem~\ref{thm:LD}. Next, define for all $n \ge 1$,
\begin{equation}\label{eq:J}
   J_n := \int_0^1 \mathbf{1}_{\{ S_n(v) \ge
   z_n\sqrt{n}\}}\, dv.
\end{equation}
Thanks to Lemma~\ref{lem:cond}, for any $u\ge 0$,
$n\ge 1$,
\begin{equation}\label{eq:project}
   \EN\left\{ J_n \, \Big|\, \F_u^n \right\}
    \geq \int_u^1 \bar\Phi\left(
    \frac{z_n\sqrt{n} -  \left( 1-\frac1n N^n_{u\to v} \right)S_n(u)}{
    \sqrt{N^n_{u\to v}\left( 2- \frac1n N^n_{u\to v}\right)}}
    \right)\,  dv.
\end{equation}
Now consider the following ``good'' events, where $n\ge 1$ is an integer,
and $\alpha\in(0,1)$ is an arbitrarily small parameter:
\begin{equation}\label{eq:good}\begin{split}
   A_{n,\alpha} &:= \left\{ \sup_{\substack{
      0\le s\le t \le 1:\\ t-s\ge\Delta_n}}
      \left| \frac{N^n_{s\to t}}{\E N^n_{s\to t}
      } - 1 \right| \le \alpha \right\},\\
   B_n(u) &:= \left\{
      S_n(u) \ge z_n\sqrt{n} \right\}.
\end{split}\end{equation}
Next is a key technical estimate.

\begin{lemma}\label{lem:tech-UB}
   Choose and fix integers $n, m\ge 1$, $u\in[0,1-\frac1m]$, and
   $\alpha\in(0,1)$. Then, $\P$-a.s.,
  \begin{equation}
      \EN\left\{ J_n \, \Big|\, \F_u^n \right\} \ge
      \frac{1}{(1+\alpha)z_n^2}\int_{0}^{z_n^2/m}
         \bar\Phi\left(\sqrt{t}\right)\, dt
         \cdot \mathbf{1}_{A_{n,\alpha}
         \cap B_n(u)}.
   \end{equation}
\end{lemma}

\begin{proof}
   Thanks to (\ref{eq:project}), for any $u\ge 0$,
   \begin{equation}\label{eq:project-good}\begin{split}
      &\EN\left\{ J_n \, \Big|\, \F_u^n \right\}\\
      &  \ge \int_{u}^1 \bar\Phi\left(
         \frac{z_n\sqrt{n} -  \left( 1-\frac1n
         N^n_{u\to v} \right)S_n(u)}{
         \sqrt{N^n_{u\to v}\left( 2- \frac1n
         N^n_{u\to v}\right)}}
         \right)\, dv \cdot \mathbf{1}_{A_{n,\alpha}
         \cap B_n(u)}.
   \end{split}\end{equation}
   We will estimate the terms inside $\bar\Phi$.
   On $B_n(u)$, we have
   \begin{equation}\begin{split}
      \frac{z_n\sqrt{n} -  \left( 1-\frac1n N^n_{u\to v} \right)S_n(u)}{
         \sqrt{N^n_{u\to v}\left( 2- \frac1n N^n_{u\to v}\right)}}
      &\le\frac{z_n\sqrt{n} -  
         \left( 1-\frac1n N^n_{u\to v}\right) z_n\sqrt{n}}{
         \sqrt{N^n_{u\to v}}}\\
      &= z_n\sqrt{ \frac{N^n_{u\to v} }{ n}}.
   \end{split}\end{equation}
   On the other hand, on $A_{n,\alpha}$,
   \begin{equation}
      N^n_{u\to v} \le (1+\alpha) n \left( 1-e^{-|v-u|}\right)
      \le (1+\alpha) (v-u) n.
   \end{equation}
   Consequently, on $A_{n,\alpha}\cap B_n(u)$,
   the preceding two displays combine to yield the following:
   \begin{equation}\begin{split}
      &\frac{z_n\sqrt{n} -  \left( 1-\frac1n N^n_{u\to v} \right)S_n(u)}{
         \sqrt{N^n_{u\to v}\left( 2- \frac1n N^n_{u\to v}\right)}}
      \le z_n\sqrt{(1+\alpha)(v-u) }.
   \end{split}\end{equation}
   Because $\bar\Phi$ is decreasing, the above
   can be plugged into (\ref{eq:project-good}) to yield:
   \begin{equation}\label{eq:cond}\begin{split}
      \EN\left\{ J_n \, \Big|\, \F_u^n \right\}
         &  \ge \int_{u}^1 \bar\Phi\left(
         z_n\sqrt{(1+\alpha)(v-u)}
         \right)\,  dv \cdot \mathbf{1}_{A_{n,\alpha}
         \cap B_n(u)}\\
      &  = \frac{1}{(1+\alpha)z_n^2}
         \int_{0}^{(1-u)(1+\alpha)z_n^2}
         \bar\Phi\left(\sqrt{t}\right)\, dt
         \cdot \mathbf{1}_{A_{n,\alpha}
         \cap B_n(u)}.
\end{split}\end{equation}
The result follows readily from this.
\end{proof}

\begin{proof}[Proof of Theorem~\ref{thm:key-UB}]
   Clearly, the following holds $\P$-a.s. on $A_{n,\alpha}$:
   \begin{equation}\begin{split}
      & \PN\left\{ {}^\exists u\in \left[ 0,\textstyle{1-\frac1m}\right]:\
         S_n(u) \ge z_n\sqrt{n} \right\}\\
      &   = \PN\left\{ \sup_{u\in\left[0,1-\frac1m\right]\cap\mathbb{Q}}
         \mathbf{1}_{A_{n,\alpha} \cap
         B_n(u)} =1 \right\}.
   \end{split}\end{equation}
   Therefore, we can appeal to Lemma~\ref{lem:tech-UB}
   to deduce that $\P$-almost surely,
   \begin{equation}\label{eq:refer-UB}\begin{split}
      & \mathbf{1}_{A_{\alpha,n}}\times
         \PN\left\{ {}^\exists u\in\left[0,\textstyle{1-\frac1m}\right]:\
         S_n(u) \ge z_n\sqrt{n} \right\}\\
      &  \le \PN\left\{ \sup_{u\in\left[0,1-\frac1m\right]\cap \mathbb{Q}}
         \EN\left\{ J_n \,\big|\, \F_u^n \right\} \ge
         \frac{1}{(1+\alpha)z_n^2}\int_{0}^{z_n^2/m}
         \bar\Phi\left(\sqrt{t}\right)\, dt.\right\}\\
      &   \le \frac{(1+\al) z_n^2}{\int_{0}^{z_n^2/m}
         \bar\Phi\left(\sqrt{t}\right)\, dt}
         \EN\left\{ J_n \right\} =\frac{(1+\al) z_n^2}{\int_{0}^{z_n^2/m}
         \bar\Phi\left(\sqrt{t}\right)\, dt}
         \bar\Phi(z_n).
   \end{split}\end{equation}
   The final line uses Doob's inequality  (under $\PN$), and
   the stationarity of $S_n(u)$. According to Corollary~\ref{cor:LD},
   with $\PN$-probability one, for all but finitely-many of the
   $n$'s, $\mathbf{1}_{A_{\alpha,n}}=1$.
   To finish, we note that
   \begin{equation}
   \int_{0}^{\infty}
         \bar\Phi\left(\sqrt{t}\right)\, dt=\frac12.
   \end{equation}
   Theorem 5.1 follows after letting $m \to \infty$ and $\alpha \to
   0$.
\end{proof}

\section{A Quenched Lower Bound}
\label{sec:quenched-LB}

\begin{theorem}\label{thm:key-LB}
   Suppose $\{ z_j\}_{j=1}^\infty$ is a sequence of
   real numbers that satisfies (\ref{eq:z}).
   Then, there exists a random variable
   $n_1$ such that $\P$-almost surely the following
   holds:
   \begin{equation}
      \PN\left\{ \sup_{t\in[0,1]}
      S_n(t) \ge z_n\sqrt{n} \right\} \ge
      \frac{1}{9} z_n^2 \bar\Phi(z_n),
      \qquad{}^\forall n\ge n_1.
   \end{equation}
\end{theorem}

We begin by proving Theorem~\ref{thm:key-LB}.

\begin{lemma}\label{lem:joint-first}
   There is some $\al_0 > 0$ so that for any fixed $\al < \al_0$,
   there exists a random variable $n_2$
   such that with $\P$-probability one,
   the following holds: For all $n\ge n_2$,
   \begin{equation}\begin{split}\label{eq:1/41}
      & \PN \left\{ S_n(u)\ge
         z_n\sqrt{n} ~,~ S_n(v)\ge z_n\sqrt{n} \right\}\\
      &\quad \le 2\exp\left( -\frac{z_n^2(1-\alpha)
         (v-u)}{4}\right)
         \bar\Phi(z_n),
   \end{split}\end{equation}
   for all $0\le u\le v\le 1$ such that $v-u\ge\Delta_n$,
   where $\Delta_n$ is defined in (\ref{eq:Delta_n}).
\end{lemma}

\begin{proof}
   In the course of our proof of Theorem~\ref{thm:key-UB}
   we observed that for any $\alpha\in(0,1)$,
   $\mathbf{1}_{A_{n,\alpha}}=1$ for all
   but a finite number of $n$'s. Thus, it suffices to derive
   the inequality of this lemma on the set $A_{n,\alpha}$.
   Recall that the latter event was defined in (\ref{eq:good}).

   By Lemma~\ref{lem:cond},
   \begin{equation}\begin{split}
      &\PN\left\{ S_n(v)\ge z_n\sqrt{n} ~,~
         S_n(u) \ge z_n\sqrt{n} \right\}\\
      & = \int_{z_n}^\infty \bar\Phi \left(
         \frac{z_n\sqrt{n}-x\sqrt{n}
         \left( 1-\frac1n N^n_{u\to v}
         \right)}{ \sqrt{N^n_{u\to v}
         \left[ 2- \frac1n N^n_{u\to v} \right]}}\right)\,
         \Phi(dx).
   \end{split}\end{equation}
  A computation shows that if $x \ge z_n$, then the function
  \begin{equation}\frac{z_n-x ( 1-u)}{\sqrt{ u
  (2-u)}},\end{equation}
  is increasing for $ u \in [0,1]$.
  On the other hand, on  $A_{n,\alpha}$, we have
   \begin{equation}N^n_{u\to v}\ge n(1-\alpha)(1-e^{-(v-u)})
   \ge n \frac12 (1-\alpha)(v-u);\end{equation} cf.\ (\ref{eq:expbound}). Therefore,
  \begin{equation}\begin{split}
      &\PN\left\{ S_n(v)\ge z_n\sqrt{n} ~,~
         S_n(u) \ge z_n \sqrt{n} \right \} \\
      & \le \int_{z_n}^\infty \bar\Phi \left(
         \frac{z_n-x
         \left(1-\frac12(1-\alpha)(v-u)
         \right)} {\sqrt{\frac12(1-\alpha)(v-u)
         \left[ 2- \frac12(1-\alpha)(v-u) \right] }}\right)\,
         \Phi(dx)\\
      & = \int_{z_n}^\infty \bar\Phi \left(
         \frac{\frac12 x(1-\alpha)(v-u)-
         (x-z_n)}{\sqrt{ \frac12 (1-\alpha)(v-u)
         }}\right)\,
         \Phi(dx)\\
      & := \mathrm{I}_1 +\mathrm{I}_2,
   \end{split}\end{equation}
   where $\mathrm{I}_1:=\int_{z_n}^{(1+\eta)z_n}\bar\Phi(\cdots)\, \Phi(dx)$,
   $\mathrm{I}_2:=\int_{z_n(1+\eta)}^\infty\bar\Phi(\cdots)\, \Phi(dx)$, and
   \begin{equation}\label{eq:eta}
      \eta := \frac{\gamma}{2} (1-\alpha)(v-u).
   \end{equation}
   $\gamma \in (0,1)$ is a parameter to be determined.
   For the estimation of $\mathrm{I}_1$, we note
   that if $x\in[z_n,z_n(1+\eta)]$, then
   $\frac12 x\sqrt{n}(1-\alpha)(v-u)-
   (x-z_n)\sqrt{n}\ge z_n\frac12(1-\alpha)(v-u)\,(1-\gamma)$,
   and we obtain the following:
   \begin{equation}\begin{split}
      \mathrm{I}_1 &  \le \int_{z_n}^\infty \bar\Phi\left(
         z_n(1-\gamma)\sqrt{\frac12(1-\alpha)(v-u)}\right)\,
         \Phi(dx)\\
      & \le \exp\left( - \frac{z_n^2(1-\gamma)^2
       (1-\alpha) (v-u)}{4}
         \right) \bar\Phi(z_n),
   \end{split}\end{equation}
   where the last line follows from (\ref{eq:Phiasymp}).
   The integral $\mathrm{I}_2$ is also easily
   estimated: Since $\bar\Phi(t)\le 1$, we have
   \begin{equation}
      \mathrm{I}_2  \le \bar\Phi\left( z_n(1+\eta)\right)
      \le \exp\left( - \eta z_n^2\right)
      \bar\Phi(z_n)  \le e^{- \eta z_n^2}
      \bar\Phi(z_n).
   \end{equation}
   We have appealed to Lemma~\ref{lem:phi-bar}
   in the penultimate inequality. Now replace $\eta$ by
   its value defined in (\ref{eq:eta})
   in order to obtain 
   \begin{equation}
     \mathrm{I}_2 \leq \exp\left( -  z_n^2 \frac{\gamma}{2}
     (1-\al)(v-u)\right) \bar\Phi(z_n)\,.
   \end{equation}
   Taking $\gamma$ to be the solution of $  \gamma=
   \frac{(1-\gamma)^2}{2}$ in $[0,1]$ we have that
   \begin{equation}
     \mathrm{I}_1 + \mathrm{I}_2 \leq 2 \exp\left( - \left(2
     -\sqrt{3}\right)
     (1 - \al)(v - u) \right) \bar\Phi(z_n) \,,
   \end{equation}
    the  result follows from the fact that $(2
     -\sqrt{3}) \leq \frac14$.
\end{proof}

\begin{proof}[Proof of Theorem~\ref{thm:key-LB}]
   We recall (\ref{eq:J}) and appeal to
   Lemma~\ref{lem:joint-first} to see that
   $\P$-a.s., for all $n\ge n_3$,
   \begin{equation}\begin{split}
      \EN \left\{ J_n^2 \right\}
         &= 2\int_0^1\int_u^1
         \PN\left\{ S_n(v)\ge z_n\sqrt{n} ~,~
         S_n(u) \ge z_n\sqrt{n} \right\}\, dv\, du\\
      & \le 2\bar\Phi(z_n) \int_0^{1-\Delta_n}\int_{u+\Delta_n}^1
         \exp\left( - \frac{z_n^2(1-\alpha) (v-u)}{4}
         \right)\, dv\, du\\
      &\quad + 2\Delta_n \bar\Phi(z_n)\\
      & \le z_n^{-2}\bar\Phi(z_n)\left[
         \frac{8}{(1-\alpha)} +\frac{2}{16}\right].
   \end{split}\end{equation}
   We have used the definition (\ref{eq:Delta_n})
   of $\Delta_n$ in the last line. Let us choose
   $\alpha$ small enough so that $8/(1-\al) + 1/8 < 9$. Then, we obtain:
   \begin{equation}
      \EN \left\{ J_n^2 \right\} \le 9 z_n^{-2}
      \bar\Phi(z_n),\text{ a.s. on }A_{\alpha,n}.
   \end{equation}
   Thus, by the Paley--Zygmund inequality,
   almost surely on $A_{\alpha,n}$,
   \begin{equation}\label{eq:key-lower}\begin{split}
      \PN \left\{ J_n>0 \right\}
      \ge \frac{\left( \EN J_n \right)^2}{ \EN J_n^2 }
      \ge \frac{1}{9}  z_n^2 \bar\Phi(z_n).
   \end{split}\end{equation}
   The theorem follows readily from this and the obvious fact that
   $\{J_n(z_n)>0\}\subseteq\{
   {}^\exists u\le 1: S_n(u)\ge z_n\sqrt{n}\}$.
\end{proof}

\section{Proof of Theorem~\ref{thm:key-est}}
\label{sec:key}

We start by proving the simpler lower bound. Fix $\alpha\in(0,1)$,
let $W_n$ denote the $\PN$-probability that
$\sup_{t\in[0,1]}S_n(t)\ge z_n\sqrt{n}$, and define $f_n :=
z_n^2\bar\Phi(z_n)$. [We will use this notation throughout the
proof.] Then, according to (\ref{eq:key-lower}), $9 W_n\ge f_n$,
$\P$-almost surely on $A_{\alpha,n}$. Theorem~\ref{thm:LD} implies
that $\P(A_{\alpha,n}^\mathsf{C})\to 1$, as $n\to\infty$. In
particular, as $n\to \infty$, $\P\{ 9 W_n\ge f_n\}=1+o(1)$. This,
and Chebyshev's inequality, together imply that $9 \E W_n
\ge(1+o(1))f_n$, which is the desired lower bound in scrambled
form. We now prove the corresponding probability upper bound of
Theorem~\ref{thm:key-est}.

Let $\Pi_n$ denote the total number of replacements to
the incremental processes $\{X_k(\cdot)\}_{k=1}^n$
during the time-interval $[0,1]$. That is,
\begin{equation}\label{eq:Pi}
   \Pi_n := \sum_{s\in(0,1]} \Delta\Pi_n(s),\
   \quad\text{where}\ \Delta\Pi_n(s) := \sum_{k=1}^n {\mathbf 1}_{\{X_k(s) -
   X_k(s-) \neq 0 \}}
\end{equation}
Because $\Pi_n$ is a Poisson random variable with mean $n$,
$\E\{ e^{t\Pi_n}\} = \exp( -n +e^t n)$ for all $t>0$.
This readily yields the following well-known Chernoff-type bound:
{F}or all $x>0$,
\begin{equation}
   \P\left\{ \Pi_n \ge x \right\} \le \inf_{t>0} \exp\left(-n+
   e^t n - t x \right) = \exp\left\{ -n-x\ln\left( \frac{x}{en} \right)\right\}.
\end{equation}
Consequently, by (\ref{eq:z}),
\begin{equation}\label{eq:Ccomp}
   \P\left( G_n^\mathsf{C} \right) \le
   e^{-n} = o\left(f_n\right),\ \text{where}\
   G_n  := \left\{ \Pi_n \le 3 n \right\},\qquad{}^\forall n\ge 1.
\end{equation}
A significant feature of the event $G_n$ is that
$\P$-almost surely,
\begin{equation}
   \mathbf{1}_{G_n} W_n
   \le 3n \P\left\{ S_n(0) \ge z_n\sqrt{n} \right\}
   = 3n\bar\Phi(z_n).
\end{equation}
(Indeed, if $G_n$ holds, then $W_n$ is the chance that the maximum
of at most $3n$ dependent Gaussian random walks exceeds $z_n\sqrt{n}$.)
Thus, we can write the almost sure [$\P$] bound,
\begin{equation}
   \mathbf{1}_{A_{\alpha,n}^\mathsf{C}}W_n
   \le  \mathbf{1}_{G_n^\mathsf{C}} +3n \bar\Phi(z_n)
   \mathbf{1}_{A_{\alpha,n}^\mathsf{C}}.
\end{equation}
Combined with (\ref{eq:refer-UB}) and (\ref{eq:1/41}) (for
suitable small $\al$), this yields
\begin{equation}
   W_n \le (2+o(1)) f_n
   + \mathbf{1}_{G_n^{\mathsf{C}}} + 3n \bar\Phi(z_n)
   \mathbf{1}_{A_{\alpha,n}^\mathsf{C}}.
\end{equation}
In this formula, $o(1)$ denotes a nonrandom term that 
goes to zero as $n$ tends to infinity.
We take expectations and appeal to Theorem~\ref{thm:LD} with
$\Delta_n:=(16z_n^2)^{-1}$ (cf.~\ref{eq:Delta_n}), as well as
(\ref{eq:Ccomp}), to deduce the following:
\begin{equation}
   \E\left\{ W_n\right\} \le (2+o(1)\,) f_n +
   \frac{8192}{\alpha^2} n z_n^2 f_n
   \exp\left( -\frac{3\alpha^3 n}{36864 z_n^2} \right).
\end{equation}
Condition (\ref{eq:z}) guarantees that the right-hand side is
asymptotically equal to $(2+o(1))f_n $, as $n\to\infty$. This
proves the theorem.\hfill$\square$

\section{Proof of Theorem~\ref{thm:erdos}}
\label{sec:erdos}

Throughout, $\log(x):=\log x := \ln(e\vee x)$,
and consider the \emph{Erd\H{o}s sequence}:
\begin{equation}\label{eq:ee}
   \ee_n:=\ee(n) :=\left\lfloor \exp\left( \frac{n}{%
   \log (n)} \right) \right\rfloor,\quad{}^\forall n\ge 1.
\end{equation}
Note that the sequence $\{\ee_j\}_{j=1}^\infty$ satisfies the following
\emph{gap property}:
\begin{equation}\label{eq:gap}
   \ee_{n+1}-\ee_n = \frac{\ee_n}{\log(n)}(1+o(1))
   =\frac{\ee_n}{\log\log(\ee_n)}(1+o(1)),\quad
   (n\to\infty).
\end{equation}
[This was noted in \ocite{erdos}*{eq.\ (0.11)}]
{F}urthermore, we can combine the truncation argument of
Erd\H{o}s \ycite{erdos} [eq.'s (1.2) and (3.4)]
with our equation (\ref{eq:LIL}) to
deduce the following: Without loss of generality,
\begin{equation}\label{eq:WLOG}
   \sqrt{\log\log(t)} \le H(t)
   \le 2\sqrt{\log\log(t)}\quad{}^\forall t\ge 1.
\end{equation}

The following is a standard consequence.

\begin{lemma}\label{lem:sumint}
   If $H$ is a nonnegative nondecreasing measurable function
   that satisfies (\ref{eq:WLOG}), then
   \begin{equation}\label{eq:sumint1}
      \EuScript{J}(H)<+\infty\ \Longleftrightarrow\
      \sum_n H^2(\ee_n) \bar\Phi(H(\ee_n))<+\infty,
   \end{equation}
   where $\EuScript{J}(H)$ is defined in (\ref{eq:JJ}).
\end{lemma}

We are ready to prove (the easier) part (i)
of Theorem~\ref{thm:erdos}.

\begin{proof}[Proof of Theorem~\ref{thm:erdos} (First Half)]
   In the first portion of our proof, we assume
   that $\EuScript{J}(H)<+\infty$, and recall that
   without loss of generality, (\ref{eq:WLOG}) is assumed
   to hold.

   It is easy to see that $\{ X_j\}_{j=1}^\infty$ are i.i.d. elements of
   $D([0,1])$---the space of \cadlag real
   paths on $[0,1]$---which implies that $n\mapsto S_n$ is a symmetric
   random walk on $D([0,1])$. In particular,
   an infinite-dimensional reflection argument
   implies that for all $n\ge 1$ and $\l>0$,
   \begin{equation}
      \P\left\{ \max_{\substack{{}\\ 1\le k\le n}}
      \sup_{t\in[0,1]} S_k(t) \ge \l \right\}
      \le 2 \P\left\{ \sup_{t\in[0,1]}
      S_n(t) \ge \l \right\}.
   \end{equation}
   See 
   \ocite{khosh}*{Lemma 3.5}
   for the details of this argument.
   Consequently, as $n\to\infty$,
   \begin{equation}\begin{split}
      &\P\left\{ \max_{\substack{{}\\ 1\le k\le \ee(n+1)}}
        \sup_{t\in[0,1]} S_k(t) \ge H(\ee_n) \sqrt{\ee_n}\right\}\\
      &\quad\le 2\P\left\{ \sup_{t\in[0,1]} S_{\ee(n+1)}(t)
         \ge H(\ee_n)\sqrt{\ee_n}\right\}\\
      &\quad\le 2\P\left\{ \sup_{t\in[0,1]} S_{\ee(n+1)}(t)
         \ge H(\ee_n) \sqrt{\ee_{n+1}}
         \left[ 1 -\frac{2+o(1)}{H(\ee_n)}
         \right]\right\}.
   \end{split}\end{equation}
   We have appealed to (\ref{eq:gap}) in the last line. At this point,
   (\ref{eq:WLOG}) and Theorem~\ref{thm:key-est} together imply that
   as $n\to\infty$,
   \begin{equation}\begin{split}
      &\P\left\{
      \max_{\substack{{}\\ 1\le k\le \ee(n+1)}}
      \sup_{t\in[0,1]} S_k(t) \ge
      H(\ee_n)\sqrt{\ee_n} \right\}\\
      &\quad\le [4+o(1)]
         H^2(\ee_n)\bar\Phi\left(H(\ee_n)
         \left[ 1 - \frac{2+o(1)}{H(\ee_n)}
         \right]\right)\\
      &\quad\le (e^4 4+o(1))
         H^2(\ee_n)\bar\Phi\left(H(\ee_n)\right),
   \end{split}\end{equation}
   the last line follows from Lemma~\ref{eq:phi-bari}.
   Lemma~\ref{lem:sumint} and the finiteness assumption
   on $\EuScript{J}(H)$ together yield the summability of
   the left-most probability in the preceding display.
   By the Borel--Cantelli lemma, almost surely for all
   but a finite number of $n$'s,
   \begin{equation}
      \max_{\substack{{}\\ 1\le k\le \ee(n+1)}}
      \sup_{t\in[0,1]} S_{k}(t) < H(\ee_n)\sqrt{\ee_n}.
   \end{equation}
   Now any $m$ can be sandwiched between $\ee_n$ and $\ee_{n+1}$ for
   some $n:=n(m)$. Hence, a.s. for all but a finite number of $m$'s,
   \begin{equation}
      \sup_{t\in[0,1]} S_m(t) \le \max_{\substack{{}\\ 1\le k\le \ee(n+1)}}
      \sup_{t\in[0,1]} S_{k}(t)<
      H(\ee_n)\sqrt{\ee_n}\le H(m)\sqrt{m}.
   \end{equation}
   This completes our proof of part (i).
\end{proof}

The remainder of this section is concerned with proving the
more difficult second part of Theorem~\ref{thm:erdos}.
We will continue to use the Erd\H{o}s sequence $\{ \ee_j\}_{j=1}^\infty$
as defined in (\ref{eq:ee}). We will also assume---still without loss
of generality---that (\ref{eq:WLOG}) holds, although now
$\EuScript{J}(H)=+\infty$.

We introduce the following notation in order to simplify the exposition:
\begin{equation}\begin{split}
   S^*_n & := \sup_{t\in[0,1]} S_{\ee(n)}(t)\\
   H_n  &:= H(\ee_n)\\
   \mathcal{I}_n & :=
      \left[ H_n\sqrt{\ee_n}, \left(H_n+\frac{14}{H_n}\right)
       \sqrt{\ee_n}\right]\\
   L_n  &:= \sum_{j=1}^n \mathbf{1}_{\{S^*_j\in\mathcal{I}_j\}}\\
   f(z)  &:= z^2 \bar\Phi(z),\quad{}^\forall z>0.
\end{split}\end{equation}

Here is a little localization lemma that states that
$\mathcal{I}_n$ and $[H_n\sqrt{\ee_n},+\infty]$ have, more or less,
the same dynamical-walk-measure.

\begin{lemma}\label{lem:I_n}
   As $n\to\infty$,
   \begin{equation}
      \left( 10^{-2}+o(1)\right) \le
      \frac{\P\left\{ S^*_n \in\mathcal{I}_n \right\}}{%
      \P\left\{ S^*_n \ge H_n\sqrt{\ee_n} \right\}}\le 1.
   \end{equation}
\end{lemma}

\begin{proof}
   Because $9^{-1}\ge 0.1$,
   Theorem~\ref{thm:key-est} implies that as $n\to\infty$,
   \begin{equation}\begin{split}
      \P\left\{ S^*_n \in\mathcal{I}_n \right\}
         & \ge \left( 0.1+o(1)\right) f(H_n)
         - \left( 2+o(1) \right) H_n^2
         \bar\Phi\left( H_n +\frac{14}{H_n} \right)\\
      & \ge \left( 0.1+o(1)\right) f(H_n)
         - \left( 2+o(1) \right) e^{-14} f(H_n).
   \end{split}\end{equation}
   (The second line holds because of Lemma~\ref{lem:phi-bar}.)
   Since $0.1-2 e^{-14}\le 0.09$, the lemma follows Theorem~\ref{thm:key-est}
   and a few lines of arithmetic.
\end{proof}

Since we are assuming that $\EuScript{J}(H)=+\infty$,
Lemmas~\ref{lem:sumint} and~\ref{lem:I_n} together
imply that as $n\to\infty$,
$\E L_n \to+\infty$. We intend to show that
\begin{equation}\label{eq:goal-ii}
   \limsup_{n\to\infty} \frac{
   \E\left\{ L_n^2\right\}}{%
   \left( \E L_n \right)^2}<+\infty.
\end{equation}
If so, then the Chebyshev inequality shows that
$\limsup_{n\to\infty} L_n/\E L_n > 0$ with positive probability.
This implies that with positive probability, $L_\infty=+\infty$,
so that the following would then conclude the proof.

\begin{lemma}
   If $\rho:=\P\{ L_\infty=+\infty\}>0$,
   then $\rho=1$, and part (ii) of Theorem~\ref{thm:erdos} holds.
\end{lemma}
\begin{proof}
   We have already observed that $n\mapsto S_n$ is a
   random walk in $D([0,1])$. Therefore, by the
   Hewitt--Savage 0--1 law, $L_\infty=+\infty$, a.s.

   Now consider
   \begin{equation}
      \mathcal{W}_n := \left\{ t\ge 0: S_{\ee(n)}(t) \vee
      S_{\ee(n)}(t-)
      \ge H_n\sqrt{\ee_n} \right\},\qquad{}^\forall n\ge 1.
   \end{equation}
   This is a random open set, and
   \begin{equation}
      \left\{ L_\infty=+\infty \right\} \subseteq
      \bigcap_{n=1}^\infty \bigcup_{m=n}^\infty
      \left\{ \mathcal{W}_m \cap [0,1] \neq\varnothing \right\}.
   \end{equation}
   More generally still, for any $0\le a<b$,
   \begin{equation}\label{eq:chi-a,b}
      \left\{ L_\infty(a,b)=+\infty \right\} \subseteq
      \bigcap_{n=1}^\infty \bigcup_{m=n}^\infty
      \left\{ \mathcal{W}_m \cap [a,b] \neq\varnothing \right\},
   \end{equation}
   where $L_n(a,b) := \sum_{j=1}^n \mathbf{1}_{\{
   \sup_{t\in [a,b]}S_{\ee(j)}(t)\in\mathcal{I}_j\}}$.
   But by the stationarity of the $\R^\infty$-valued
   process $t\mapsto S_\bullet(t)$, $L_\infty(a,b)$ has
   the same distribution as $L_\infty(0,b-a)$, and this means
   that with probability one,
   $L_\infty(a,b)=+\infty$ for all rational $0\le a<b$.
   Therefore, according to (\ref{eq:chi-a,b}),
   \begin{equation}
      \P\left\{ \bigcap_{n=1}^\infty \bigcup_{m=n}^\infty
      \left\{ \mathcal{W}_m \cap [a,b] \neq\varnothing
      \right\}\right\}=1.
   \end{equation}
   This development shows that for any $n$,
   $\mathcal{W}^n:=\cup_{m\ge n}\mathcal{W}_m$ is a random open set
   that is a.s. everywhere dense. Thanks to the Baire
   category theorem, $\mathcal{W}:=\cap_n \mathcal{W}^n \, \cap \, [0,1]$
   is [a.s.] uncountable.
   Now any $t\in\mathcal{W}\cap[0,1]$
   satisfies the following:
   \begin{equation}
      S_\ell(t)\vee S_\ell(t-) \ge H(\ell)\sqrt{\ell},\
      \text{for infinitely many $\ell$'s}.
   \end{equation}
   On the other hand, the jump structure of the Poisson clocks
   tells us that $\mathcal{J}:=\cup_{\ell\ge 1}
   \{t\ge 0:S_\ell(t)\neq S_\ell(t-)\}$ is [a.s.]
   denumerable. Because $\mathcal{W}$ is uncountable [a.s.],
   any $t\in\mathcal{W}\cap\mathcal{J}^\mathsf{C}$
   satisfies assertion (ii) of Theorem~\ref{thm:erdos}.
\end{proof}

We now begin working toward our proof of (\ref{eq:goal-ii}).
We write
\begin{align}
   \E\left\{ L_n^2 \right\} &= \E L_n
      + 2 \sum_{i=1}^{n-1} \sum_{j=i+1}^n \mathcal{P}_{i,j},
   \intertext{where}
   \mathcal{P}_{i,j} & = \P\left\{ S^*_i \in\mathcal{I}_i ~,~
      S^*_j \in \mathcal{I}_j \right\},\qquad{}^\forall
      i>j\ge 1.
\end{align}
In estimating $\mathcal{P}_{i,j}$, our first observation is the following.
\begin{lemma}\label{lem:Pij}
   There exists a finite and positive
   universal constant $K_{\ref{lem:Pij}}$ such that for
   all $j>i\ge 1$,
   \begin{align}
      \mathcal{P}_{i,j} &\le K_{\ref{lem:Pij}}
         \P\left\{ S^*_i \in \mathcal{I}_i\right\} \mathcal{Q}_{i,j}\label{eq:Pij}
      \intertext{where}
      \mathcal{Q}_{i,j} &:= f\left( H_j\sqrt{\frac{\ee_j}{\ee_j-\ee_i}}
         - H_i\sqrt{\frac{\ee_i}{\ee_j-\ee_i}} -
         \frac{14}{H_i}\sqrt{\frac{\ee_i}{\ee_j-\ee_i}} \right).
         \label{eq:Qij}
   \end{align}
\end{lemma}

\begin{proof}
   Recall that $n\mapsto S_n$ is a random walk on
   $D([0,1])$. Therefore,
   \begin{equation}\begin{split}
      \mathcal{P}_{i,j} &\le \P\left\{ S^*_i\in\mathcal{I}_j\right\}\\
      &\quad \times \P\left\{ \sup_{t\in[0,1]} \left(
         S_{\ee_j}(t) - S_{\ee_i}(t) \right)
         \ge H_j\sqrt{\ee_j} - \sqrt{\ee_i}\left[
         H_i+ \frac{14}{H_i}  \right] \right\}\\
      & = \P\left\{ S^*_i\in \mathcal{I}_i \right\}
         \P\left\{ \sup_{t\in[0,1]}
         S_{\ee_j-\ee_i}(t) \ge H_j\sqrt{\ee_j} - \sqrt{\ee_i}\left[
         H_i+ \frac{14}{H_i}  \right] \right\}.
   \end{split}\end{equation}
   Therefore, Theorem~\ref{thm:key-est} will do the rest, once
   we check that uniformly for all $j>i$,
   \begin{equation}\label{eq:WLOG-CHECK}
      \frac{H_j\sqrt{\ee_j}}{\sqrt{\ee_j-\ee_i}}
      = o\left( \sqrt{\frac{\ee_j-\ee_i}{\log(\ee_j-\ee_i)}}\right)
      \qquad(i\to\infty).
   \end{equation}
   Equivalently, we wish to prove that uniformly for all $j>i$,
   \begin{equation}
      H_j \sqrt{\ee_j} = o\left(
      \frac{\ee_j-\ee_i}{\sqrt{\log(\ee_j-\ee_i)}}
      \right)\qquad (i\to\infty).
   \end{equation}
   By (\ref{eq:WLOG}), the left-hand side is bounded above as follows:
   \begin{equation}\label{eq:WLOG-LHS}
      H_j \sqrt{\ee_j} \le (2+o(1))\sqrt{\ee_j\log\log\ee_j}=
      O\left( \sqrt{\ee_j \log j}\right),\quad(j\to\infty).
   \end{equation}
   On the other hand,
   \begin{equation}\label{eq:WLOG-RHS}
      \frac{\ee_j-\ee_i}{\sqrt{\log(\ee_j-\ee_i)}} \ge
      \frac{\ee_j-\ee_i}{\sqrt{\log \ee_j}}
      = (\ee_j-\ee_i) \sqrt{\frac{\log j}{j}}.
   \end{equation}
   In light of (\ref{eq:WLOG-LHS}) and (\ref{eq:WLOG-RHS}),
   (\ref{eq:WLOG-CHECK})---and hence the lemma---is proved
   once we verify that as $i\to\infty$, $\sqrt{j\ee_j}=o(\ee_j-\ee_i)$
   uniformly for all $j>i$. But this follows from the gap condition
   of the sequence $\ee_1,\ee_2,\ldots$. Indeed, (\ref{eq:gap})
   implies that uniformly for all $j>i$,
   \begin{equation}\label{eq:ej-ej-diff}
      \ee_j - \ee_i  \ge \ee_j -\ee_{j-1} =
      (1+o(1)) \frac{\ee_j}{\log j} \quad(i\to\infty).
   \end{equation}
   So it suffices to check that as $j\to\infty$,
   $\sqrt{j\ee_j} = o(\ee_j/\log j)$, which is a trivial matter.
\end{proof}
Motivated by the ideas of P\'al Erd\H{o}s \ycite{erdos},
we consider the
size of $\mathcal{Q}_{i,j}$
on three different scales, where $\mathcal{Q}_{i,j}$ is defined in
(\ref{eq:Qij}). The mentioned scales are based on the size of the ``correlation
gap,'' $(j-i)$.  Our next three lemmas reflect this viewpoint.

\begin{lemma}\label{lem:P1}
   There exists a finite and positive universal constant $K_{\ref{lem:P1}}$
   such that for all integers $i$ and $j>i+\left[\log i\right]^{10}$,
   \begin{equation}
      \mathcal{Q}_{i,j} \le K_{\ref{lem:P1}}
      \P\left\{ S^*_j \in\mathcal{I}_j \right\}.
   \end{equation}
\end{lemma}
\begin{proof}
   We will require the following consequence
   of (\ref{eq:gap}): Uniformly for all integers $j>i$,
   \begin{equation}\label{eq:ej-ei-diff}
      \ee_j - \ee_i = \sum_{l=i}^{j-1} \left(
      \ee_{l+1}-\ee_l \right) \ge
      \frac{ (j-i) \ee_i}{\log i}(1+o(1))
      \qquad(i\to\infty).
   \end{equation}
   Now we proceed with the proof.

   Since $\ee_j/(\ee_j-\ee_i)\ge 1$, (\ref{eq:Qij}) implies that
   \begin{equation}
      \mathcal{Q}_{i,j} \le f\left( H_j
      - \sqrt{\frac{\ee_i}{\ee_j-\ee_i}}\left[H_i+
      \frac{14}{H_i}\right]\right).
   \end{equation}

   We intend to prove that uniformly for
   every integer $j\ge i+\left[\log i\right]^{10}$,
   \begin{equation}\label{eq:goal-P1}
      \sqrt{\frac{\ee_i}{\ee_j-\ee_i}}\left[H_i+
      \frac{14}{H_i}\right]
      = O\left( H_j^{-1} \right)\qquad(i\to\infty).
   \end{equation}
   Given this for the time being, we finish the proof as
   follows: Note that the preceding display
   and (\ref{eq:phi-bari}) together prove that uniformly for
   every integer $j\ge i+\left[\log i\right]^{10}$,
   $\mathcal{Q}_{i,j} = O(f(H_j))$ as $i\to\infty$. According to
   Theorem~\ref{thm:key-est}, for this range of $(i,j)$,
   $\mathcal{Q}_{i,j}=O(\P\{ S^*_j\ge H_j\sqrt{\ee_j}\})$.
   Thanks to Lemma~\ref{lem:I_n},
   this is $O(\P\{S^*_j\in\mathcal{I}_j\})$.
   The result follows easily from this,
   therefore it is enough to derive (\ref{eq:goal-P1}).

   Because of (\ref{eq:WLOG}), equation (\ref{eq:goal-P1}) is
   equivalent to the following: Uniformly for
   every integer $j\ge i+\left[\log i\right]^{10}$,
   \begin{equation}
      \frac{\ee_i (\log i)(\log j)}{\ee_j-\ee_i}
      = O(1)\qquad(i\to\infty).
   \end{equation}
   But thanks to (\ref{eq:ej-ei-diff}), uniformly for all
   integers $j>i+\left[\log i\right]^{10}$, the left-hand side is
   at most
   \begin{equation}
      (1+o(1)) \frac{\left[\log i\right]^2\log\left( i+\left[\log i\right]^{10}
\right)}{
      \left[\log i\right]^{10}} =o(1)\qquad(i\to\infty).
   \end{equation}
   This completes our proof.
\end{proof}

\begin{lemma}\label{lem:P2}
   Uniformly for all integers $j\in[i+\log i,i + \left[\log i\right]^{10}]$,
   \begin{equation}
      \mathcal{Q}_{i,j} \le i^{-\frac14+o(1)}
      \qquad(i\to\infty).
   \end{equation}
\end{lemma}

\begin{proof}
   Whenever $j>i$, we have $H_j\ge H_i$.
   Thus, the (eventual) monotonicity of $f$ implies
   that as $i\to\infty$, the following holds
   uniformly for all $j>i$:
   \begin{equation}\begin{split}
      \mathcal{Q}_{i,j} &\le f\left( H_i \left[
         \sqrt{\frac{\ee_j}{\ee_j-\ee_i}}-\sqrt{\frac{\ee_i}{\ee_j-\ee_i}}
         -\frac{14}{H_i^2}\sqrt{\frac{\ee_i}{\ee_j-\ee_i}}
         \right]\right)\\
      & = f\left( H_i \left[
         \frac{\sqrt{\ee_j-\ee_i}}{\sqrt{\ee_j}+\sqrt{\ee_i}}
         -\frac{14}{H_i^2}\sqrt{\frac{\ee_i}{\ee_j-\ee_i}}
         \right]\right)\\
      &\le f\left( H_i \left[
         \frac{\sqrt{\ee_j-\ee_i}}{\sqrt{\ee_j}+\sqrt{\ee_i}}
         -\frac{14+o(1)}{H_i^2}\sqrt{\frac{\ee_i\log j}{\ee_j}}
         \right]\right).
   \end{split}\end{equation}
   [The last line relies on (\ref{eq:ej-ej-diff}).] According to
   (\ref{eq:WLOG}), and after appealing to the trivial inequality
   that $\ee_j\ge\ee_i$, we arrive at the following:
   As $i\to\infty$, then uniformly for all integers
   $j\in[i+\log i,i + \left[\log i\right]^{10}]$,
   \begin{equation}\label{eq:Qij-LB}\begin{split}
      \mathcal{Q}_{i,j} & \le f\left( \frac{1+o(1)}{2}\sqrt{\log i} \left[
         \sqrt{\frac{\ee_j-\ee_i}{\ee_j}}
         -O\left(\frac{\sqrt{\log j}}{\log i}\right)
         \right]\right)\\
      & \le f\left( \frac{1+o(1)}{2}\left[ \sqrt{\log i}
         \sqrt{\frac{\ee_j-\ee_i}{\ee_j}}
         -O(1)\right] \right)\\
      & \le \exp\left\{ -\frac{1+o(1)}{4}\left[
         \frac{\ee_j-\ee_i}{\ee_j} \right]  \log i  \right\}.
   \end{split}\end{equation}
   [The last line holds because of the first inequality
   in (\ref{eq:Phiasymp}).] On the other hand, uniformly for all
   $j\ge i+\log i$,
   \begin{equation}\begin{split}
      \frac{\ee_j}{\ee_i} & =\exp\left( \frac{j}{\log j}
         -\frac{i}{\log i} \right) \\
      & \ge \exp\left( \frac{i+\log i}{\log\left( i+\log i\right)}
         -\frac{i}{\log i} \right) \\
      & \ge 2+o(1) \qquad(i\to\infty).
   \end{split}\end{equation}
   Consequently, $\ee_j-\ee_i\ge (1+o(1))\ee_j$.
   This and (\ref{eq:Qij-LB}) together yield the lemma.
\end{proof}

\begin{lemma}\label{lem:P3}
   Uniformly for all integers $j\in(i,i+\log i]$,
   \begin{equation}
      \mathcal{Q}_{i,j} \le \exp\left\{ -\frac{1+o(1)}{4e} (j-i) \right\}
      \qquad(i\to\infty).
   \end{equation}
\end{lemma}
\begin{proof}
   Equation (\ref{eq:ej-ei-diff})
   tell us that uniformly for all integers $j> i$, and as $i\to\infty$, $\ee_j-\ee_i \ge
   (1+o(1)) \ee_i(j-i)/\log i$.
   On the other hand, for  $j\in(i,i+\log i]$,
   \begin{equation}
      \frac{\ee_j}{\ee_i} = \exp\left( \frac{j}{\log j}
      -\frac{i}{\log i }\right)
      \le \exp\left( \frac{j-i}{\log i} \right)\le e.
   \end{equation}
   The preceding two displays together yield that uniformly
   for all integers $j\in(i,i+\log i]$,
   $\ee_j^{-1}(\ee_j-\ee_i) \ge (1+o(1))(j-i)/(e\log i)$
   ($i\to\infty$).
   The lemma follows from this and (\ref{eq:Qij-LB}).
\end{proof}

We are ready to commence with the following.
\begin{proof}[Proof of Theorem~\ref{thm:erdos}]
   Recall that $\E L_n \to\infty$, and our
   goal is to verify (\ref{eq:goal-ii}). According to
   Lemma~\ref{lem:Pij}, given any two positive integers $n>k$,
   \begin{equation}\begin{split}
      \E\left\{ \left(L_n-L_k\right)^2 \right\} & =
         \E\left\{L_n -L_k \right\}
         + 2 \sum_{i=k}^{n-1} \sum_{j=i+1}^n \mathcal{P}_{i,j}\\
      & \le \E L_n +2 K_{\ref{lem:Pij}}
         \sum_{i=k}^{n-1} \sum_{j=i+1}^n
         \P\left\{ S^*_i \in\mathcal{I}_i\right\}
         \mathcal{Q}_{i,j}.
   \end{split}\end{equation}
   We split the double-sum according to whether
   $j>i+\left[\log i\right]^{10}$, $j\in(i+\log i,i+\left[\log i\right]^{10}]$,
   or $j\in(i,i+\log i]$ and respectively
   apply Lemmas~\ref{lem:P1},~\ref{lem:P2}, and
   \ref{lem:P3} to deduce the existence of
   an integer $\nu\ge 1$ such that for all $n>\nu$,
   \begin{align}
      &\E\left\{ \left(L_n -L_{\nu}\right)^2\right\} \nonumber\\
      \begin{split}
         &\quad \le \E L_n +2 K_{\ref{lem:Pij}}K_{\ref{lem:P1}}
            \mathop{\sum\sum}\limits_{%
            \substack{ \nu\le i\le n \\
            n\ge j>i+\left[\log i\right]^{10}}}\P\left\{
            S^*_i\in\mathcal{I}_i\right\}
            \P\left\{ S^*_j\in\mathcal{I}_j\right\}\\
         &\qquad + 2 K_{\ref{lem:Pij}} \mathop{\sum\sum}\limits_{%
            \substack{ \nu\le i\le n \\
            j\in \left(i+\log i,i+\left[\log i\right]^{10}\right]}}
            i^{-1/8} \P\left\{ S^*_i\in\mathcal{I}_i \right\}\\
         &\qquad + 2 K_{\ref{lem:Pij}} \mathop{\sum\sum}\limits_{%
            \substack{ \nu\le i\le n \\
            j\in \left(i,i+\log i\right]}}
            e^{ -(j-i)/12}
            \P\left\{ S^*_i\in\mathcal{I}_i \right\}.
   \end{split}\end{align}
   Since $\E L_n\to\infty$, the above is at most
   $2 K_{\ref{lem:Pij}}K_{\ref{lem:P1}}(1+o(1))(\E L_n)^2$
   as $n\to\infty$. This proves our claim (\ref{eq:goal-ii}).
\end{proof}

\end{document}